\documentclass[12pt,twoside]{article}
\usepackage{amsmath}
\usepackage{amsfonts}
\usepackage{amsthm}
\usepackage{amssymb}
\usepackage{graphicx}
\usepackage{multicol}

\begin{document}
\title{{\normalsize{\bf Whitehead aspherical conjecture via ribbon sphere-link}}}
\author{{\footnotesize Akio KAWAUCHI}\\
{\footnotesize{\it Osaka Central Advanced Mathematical Institute, Osaka Metropolitan University }}\\
{\footnotesize{\it Sugimoto, Sumiyoshi-ku, Osaka 558-8585, Japan}}\\
{\footnotesize{\it kawauchi@omu.ac.jp}}}
\date\, 
\maketitle
\vspace{0.25in}
\baselineskip=15pt
\newtheorem{Theorem}{Theorem}[section]
\newtheorem{Conjecture}[Theorem]{Conjecture}
\newtheorem{Lemma}[Theorem]{Lemma}
\newtheorem{Sublemma}[Theorem]{Sublemma}
\newtheorem{Proposition}[Theorem]{Proposition}
\newtheorem{Corollary}[Theorem]{Corollary}
\newtheorem{Claim}[Theorem]{Claim}
\newtheorem{Definition}[Theorem]{Definition}
\newtheorem{Example}[Theorem]{Example}

\begin{abstract} 
Whitehead aspherical conjecture says that every connected subcomplex of every aspherical 
2-complex is aspherical. 
By an argument on ribbon sphere-links, it is confirmed that the conjecture is true 
for every contractible finite 2-complex. In this paper, by generalizing this argument, this conjecture is confirmed to be true for every aspherical 
2-complex.

\phantom{x}

\noindent{\footnotesize{\it Keywords: infinite ribbon sphere-link,\, aspherical 2-complex,\, 
contractible 2-complex,\, Whitehead aspherical conjecture}} 

\noindent{\footnotesize{\it Mathematics Subject classification 2010}:57Q45,\, 57M20 }
\end{abstract}

\maketitle

\phantom{x}

\noindent{\bf 1. Introduction} 

A {\it finite} or {\it infinite}  2-complex is a finite or countably-infinite CW 2-complex 
constructed from a connected finite or countably-infinite graph by attaching a finite or 
countably-infinite system of 2-cells with attaching maps, respectively. 
A  2-complex is homotopy equivalent to a simplicial  2-complex 
constructed from a simplicial  graph by attaching 2-cells with 
simplicial approximations of the attaching maps. By this homotopy equivalence, 
every subcomplex of a 2-complex is also homotopy equivalent to a simplicial 
subcomplex of the simplicial 2-complex (cf. Spanier \cite{Spanier} for a general reference). 
A path-connected space $X$ is {\it aspherical} if the universal cover $\tilde X$ of $X$ is contractible (i.e., homotopy equivalent to a point). In particular, a connected 2-complex $P$ is aspherical if and only if the second homotopy group $\pi_2(P,v)=0$. 
The Whitehead asphericity conjecture is the following conjecture (see (\cite{Howie, White}).

\phantom{x}

\noindent{\bf Conjecture~1.}
Every connected subcomplex of any aspherical 2-complex is aspherical. 

\phantom{x}

The purpose of this paper is to show that Conjecture~1 is yes. That is, 

\phantom{x}

\noindent{\bf Theorem~1.1.} Whitehead Aspherical Conjecture is true. 

\phantom{x}

A 2-complex $P$ is {\it locally finite} if every 1-cell of $P$ attaches only to a finite number of 2-cells of $P$. Conjecture 1 reduces to the following conjecture for every contractible locally  finite 2-complex.

\phantom{x}

\noindent{\bf Conjecture~2.}
Every connected subcomplex of every contractible locally finite 2-complex is aspherical. 

\phantom{x} 

In Section~2, the claim of Conjecture~2 $\Rightarrow$ Conjecture~1 is shown. 
In \cite{KA}, Conjecture~2 for every contractible finite 2-complex is 
confirmed. 
In this paper, the argument for an infinite 2-complex becomes the main argument. 
The {\it 2-complex} of a group presentation 
\[GP=<x_1, x_2,  \dots , x_n,  \dots  |\, r_1, r_2,  \dots , r_m,  \dots >\] 
is the connected 2-complex constructed from a graph with fundamental group isomorphic to 
the free group $<x_1, x_2,  \dots , x_n,  \dots >$ on the generators $x_i\,(i=1,2, \dots ,n, \dots )$ 
by attaching 2-cells with attaching maps given by the relators $r_j\,(j=1,2, \dots ,m, \dots )$, 
where note that this 2-complex is a connected graph for the empty relator.  
Up to cellular homotopy equivalences, every connected 2-complex $P$ and the connected subcomplexes of $P$ can be uniquely considered as 
the 2-complex and the subcomplexes of a group presentation $GP$, where a subcomplex of $GP$ is the 2-complex of the group presentation of a sub-presentation 
\[<x_{i_1}, x_{i_2},  \dots , x_{i_s},  \dots  |\, r_{j_1}, r_{j_2},  \dots , r_{j_t},  \dots >.\] 
A group presentation $GP$ is {\it locally finite} if every generator $x_i$ appears only in a finite number of the relators $r_j\,(j=1,2, \dots ,m, \dots )$. The 2-complex of a locally finite group presentation $GP$ can be taken as a connected locally finite  2-complex. 
A group presentation $GP$ is 
a {\it homology-trivial unit-group presentation} 
if $GP$ is a presentation of the unit group $\{1\}$ and the relator word $r_j$ is equal to the generator $x_j$ for every $j$ in the abelianized free abelian group 
$<x_1, x_2,  \dots  x_n,  \dots >^a$  of the free group 
$<x_1, x_2,  \dots  x_n,  \dots >$ with  as basis the generators $x_i\,(i=1,2, \dots ,n, \dots )$ of $GP$. 
Note that the 2-complex of a homology-trivial unit-group presentation is always contractible.

In Section~3, Conjecture~2 reduces to the following conjecture for the 2-complex of every 
homology-trivial unit-group presentation.

\phantom{x}

\noindent{\bf Conjecture~3.}
Every subcomplex of every homology-trivial unit-group presentation is aspherical. 

\phantom{x} 

The claim that Conjecture~3 $\Rightarrow$ Conjecture~2 is shown there. 
For this purpose, after observations on base changes of a free  group and 
a free abelian group of possibly infinite ranks, 
it is shown that if the 2-complex of a locally finite group presentation $GP$
is contractible, then there is a base change $x'_i\, (i=1,2, \dots ,n, \dots )$ in the free group 
$<x_1, x_2,  \dots , x_n,  \dots >$ with as basis the generators 
$x_i(i=1,2, \dots ,n, \dots )$ of  $GP$ so that the resulting group presentation 
\[GP'=<x'_1, x'_2,  \dots , x'_n,  \dots  |\, r'_1, r'_2,  \dots , r'_m,  \dots >\]
is a homology-trivial unit-group presentation (see Lemma~3.2). 
This means that there is a cellular-homotopy 
equivalence from every contractible locally finite 2-complex $P$ to the 2-complex $P'$ of 
a homology-trivial unit-group presentation $GP$ inducing a cellular-homotopy 
equivalence from the subcomplexes of $P$ to the subcomplexes of $P'$ (see Corollary~3.3). 

In Section~4, a (possibly infinite) sphere-link (namely, an $S^2$-link) $L$ in the 4-space ${\mathbf R}^4$ is discussed. 
The closed complement of $L$ in ${\mathbf R}^4$ is denoted by $E(L)$.
It is shown there that for every homology-trivial unit-group presentation 
$GP$, a ribbon $S^2$-link $L$ in ${\mathbf R}^4$ is constructed so that 
the fundamental group $\pi_1(E(L),v)$ is isomorphic to the 
free group $<x_1, x_2,  \dots , x_n,  \dots >$ of 
the generators $x_i\, (i=1,2, \dots ,n, \dots )$ of $GP$ by an isomorphism sending  a meridian system of $L$ in $\pi_1(E(L),v)$ to the relator word system $r_j\,(j=1,2, \dots ,m, \dots )$ of $GP$ (see Lemma~4.1). 
It is also observed there that a  ribbon $S^2$-link $L$ in ${\mathbf R}^4$  contains 
canonically a ribbon disk-link $L^D$ in the upper-half 4-space ${\mathbf H}^4$  so that 
the fundamental group $\pi_1(E(L),v)$ is canonically identified with 
 the fundamental group $\pi_1(E(L^D),v)$ for the closed exterior $E(L^D)$ 
of $L^D$ in ${\mathbf H}^4$  (see Lemma~4.2 and Corollary~4.3).

In Section~5, it is shown that $E(L^D)$ is always aspherical and 
every 1-full subcomplex $P'$ of the 2-complex $P$ of a homology-trivial unit-group presentation $GP$ is homotopy equivalent to the closed exterior $E(L^D)$ of 
a ribbon disk-link $L^D$ in ${\mathbf H}^4$, where a {\it 1-full subcomplex} $P'$ of $P$ 
is a subcomplex of $P$ containing the $1$-skelton $P^1$ of $P$. 
Then Conjecture~3 is confirmed to be true and the proof of Theorem~1.1 
is completed. 

The author mentions here that there is a preprint by Pasku \cite{Pasku} claiming the same result, which is a purely group-theoretic argument much different from the current argument.

\phantom{x}

\noindent{\bf 2. Reducing to the conjecture for a contractible locally finite 2-complex}

In this section, it is explained that Conjecture~1 (Whitehead Asphericity Conjecture) is obtained 
from the following conjecture. 

\phantom{x}

\noindent{\bf Conjecture~2.} 
Every connected subcomplex of every contractible locally finite 2-complex is aspherical.

\phantom{x}

For this reduction, the following three lemmas are used.

\phantom{x}

\noindent{\bf Lemma~2.1} If every connected finite subcomplex of a contractible 2-complex $P$ is aspherical, then every connected subcomplex of $P$ is aspherical. 

\phantom{x}

\noindent{\bf Lemma~2.2.} If every connected subcomplex of every contractible 2-complex is aspherical, then every connected subcomplex $Q$ of every aspherical 2-complex $P$ is aspherical.

\phantom{x}

\noindent{\bf Lemma~2.3.} Every connected finite subcomplex of a connected infinite
2-complex $P$ is a subcomplex of a connected locally finite 2-complex $P'$ homotopy equivalent to $P$.

\phantom{x}

Proof of Lemma~2.1 is done as follows.
.

\phantom{x}

\noindent{\bf Proof of Lemma~2.1.} Let $Q$ be any given connected subcomplex of a contractible 
2-complex $P$. Let $f:S^2 \to  |Q|$ be a map from the 2-sphere $S^2$ to the polyhedron $|Q|$. 
For a simplicial 2-complex $Q$, the topology of $|Q|$ is 
the topology coherent with the simplexes of $Q$ (see \cite[p.111]{Spanier}), so that 
the image $f(S^2)$ is in the polyhedron $|Q^f|$ of a connected finite subcomplex $Q^f$ of $Q$. 
By assumption, $Q^f$ is aspherical, so that the map $f:S^2 \to  |Q^f|$ defined by the original 
map $f$ is null-homotopic in $|Q^f|$ and hence in $|Q|$, so that $Q$ is aspherical. 
$\square$

\phantom{x}

Proof of Lemma~2.2 is done as follows.

\phantom{x}

\noindent{\bf Proof of Lemma~2.2.} 
Let $P$ be an aspherical 2-complex, and $Q$ any connected subcomplex of $P$.
Since the universal cover $\tilde P$ of $P$ is a contractible 2-complex, 
the subcomplex $Q$ lifts to a subcomplex $\tilde Q$ of the contractible 2-complex$\tilde P$. 
and any connected component $\tilde Q_1$ of the subcomplex $\tilde Q$ is aspherical 
by assumption. 
Since the second homotopy group is independent of a covering by the lifting property (cf. \cite{Spanier}), $Q$ is aspherical. 
$\square$

\phantom{x}

Proof of Lemma~2.3 is done as follows.

\phantom{x}

\noindent{\bf Proof of Lemma~2.3.} Let $P$ be a connected infinite 2-complex, and 
$P_0$ any given connected finite subcomplex of $P$. 
Let 
\[P_0\subset P_1\subset P_2\subset \dots \subset P_n\subset \dots \]
be a sequence of connected finite subcomplexes $P_i\, (i=0,1,2, \dots ,n, \dots )$ of $P$ such that $P=\cup_{i=0}^{+\infty} P_i$. 
Let $P_i=P_{i-1}\cup J_i$ for a subcomplex $J_i$ of $P_i$ 
with $\gamma_i=P_{i-1}\cap J_i$ a graph for all $i$. 
Triangulate the rectangle $a\times[0,1]$ for every 1-simplex $a$ of $\gamma_i$ by introducing a diagonal and regard the product $\gamma_i\times[0,1]$ as a 2-complex. 
To construct a desired 2-complex $P'$, make the connected finite subcomplexes $J_i\, (i=1,2,3, \dots ,n, \dots )$ disjoint. 
Let $P'_i=P_{i-1}\cup \gamma_i\times[0,1]$ be the 2-complex obtained from the subcomplexes $P_{i-1}$ and $\gamma_i\times[0,1]$ by identifying $\gamma_i$ in $P_{i-1}$ with 
$\gamma_i\times 0$ and $\gamma_i\times1$ with $\gamma_i$ in  $J_i$ in canonical ways. 
The sequence 
\[P_0=P'_0\subset P'_1\subset P'_2\subset \dots \subset P'_n\subset \dots \]
of connected finite subcomplexes $P'_i\, (i=0, 1,2, \dots ,n, \dots )$ is obtained. 
By construction, $P'=\cup_{i=0}^{\infty}P'_i$ is a connected locally finite 2-complex containing 
$P_0$ as a subcomplex and homotopy equivalent to $P$. $\square$

\phantom{x}

Conjecture~1 is obtained from Conjecture~2 as follows.

\phantom{x}

\noindent{\bf 2.4: Proof of Conjecture~2 $\mathbf \Rightarrow$ Conjecture~1.} 
By assuming Conjecture~2, it suffices to show that every connected finite subcomplex $Q$ 
of every contractible 2-complex $P$ is aspherical. Because this claim means 
by Lemma~2.1 that every connected subcomplex of every contractible 2-complex $P$ is aspherical, which also means by Lemma~2.2 
that every connected subcomplex of every aspherical 2-complex is aspherical, 
confirming Conjecture~1. If $Q$ is a connected finite subcomplex 
of a contractible 2-complex $P$, then $Q$ is a subcomplex of a contractible locally finite 
2-complex $P'$ homotopy equivalent to $P$ by Lemma~2.3, so that 
$Q$ is aspherical by Conjecture~2. This completes the proof of 
Conjecture~2 $\mathbf \Rightarrow$ Conjecture~1. $\square$

\phantom{x}

\noindent{\bf 3. Reducing to the conjecture for the 2-complex of a homology-trivial unit-group 
presentation}

In this section, it is explained that Conjecture~2 is obtained 
from the following conjecture. 

\phantom{x}

\noindent{\bf Conjecture~3.}
Every subcomplex of every homology-trivial unit-group presentation is aspherical. 

\phantom{x} 

A {\it base change} of a free group $<x_1, x_2,  \dots  x_n,  \dots >$ with 
basis $x_i\,(i=1,2, \dots ,n, \dots )$ 
is a consequence of a finite number of the following operations, called 
{\it Nielsen transformations} (see \cite{MKS}): 

\medskip

\noindent{(1)} Exchange two of $x_i\,(i=1,2, \dots ,n, \dots )$,

\noindent{(2)} Replace an $x_i$ by $x_i^{-1}$,

\noindent{(3)} Replace an $x_i$ by $x_i x_j\,(i\ne j)$. 

\phantom{x}

A {\it base change} of a free abelian group $\mathbf A$ on a basis $a_i\,(i=1,2, \dots ,n, \dots )$ is a consequence of a finite number of the following operations: 

\medskip

\noindent{(1)} Exchange two of $a_i\,(i=1,2, \dots ,n, \dots )$,

\noindent{(2)} Replace an $a_i$ by $-a_i$,

\noindent{(3)} Replace an $a_i$ by $a_i+a_j\,(i\ne j)$.

\phantom{x}

The following lemma is well-known for a finite rank free abelian group $\mathbf A$. 

\phantom{x}

\noindent{\bf Lemma~3.1.} Let $\mathbf A$ be a free abelian group with a countable basis 
$a_i\,(i=1,2, \dots ,n, \dots )$. Let $b_i\,(i=1,2, \dots ,n, \dots )$ be a countable basis of $\mathbf A$ 
such that every column vector and every row vector of the base change matrix $C$ given by 
\[(b_1b_2 \dots  b_n \dots )=(a_1a_2 \dots  a_n \dots )C\] 
have only a finite number of non-zero entries. 
Then there is a base change of $\mathbf A$ on $a_i\,(i=1,2, \dots ,n, \dots )$
such that $C$ is the block sum $(1)\oplus C'$ for a matrix $C'$. 

\phantom{x}

\noindent{\bf Proof of Lemma~3.1.} For every $j\,(j=1,2, \dots ,n, \dots )$, let 
\[b_j=c_{1j}a_1+c_{2j}a_2+ \dots  +c_{nj} a_n+ \dots \]
be a linear combination with $(i,j)$ entries $c_{ij}$ of $C$ which are $0$ except for a finite number of 
$i\,(i=1,2, \dots ,n, \dots )$. Note that for every $j$, 
the non-zero integer system of $c_{1j}, c_{2j},  \dots , c_{nj},  \dots $ is a coprime integer system. 
By a base change (1), assume that $c_{11}$ is the smallest positive integer 
in the integers $|c_{i1}|$ for all $i$. 
For $i>1$, write $c_{i1}= \tilde c_{i1} c_{11}+d_{i1}$ for $0\leq d_{i1}<c_{11}$. 
By a base change on $a_i\,(i=1,2, \dots ,n, \dots )$, assume that 
\[b_1=c_{11}a_1+d_{21}a_2+ \dots + d_{n1} a_n+ \dots .\]
By continuing this process, it can be assumed that $b_1=a_1$.
Note that there is a positive integer $m\geq 2$ such that 
$c_{1 j}=0$ for all $j$ with $j>m$. 
Consider the linear combination 
\[b_2=c_{12}a_1+c_{22}a_2+ \dots + c_{n2} a_n+ \dots .\]
Note that the non-zero integer system of $c_{22}, c_{32}, \dots , c_{n2},  \dots $ is coprime. 
Otherwise, there is a prime common divisor $p>1$, so that $b_1$ and $b_2$ would be 
${\mathbf Z}_p$-linearly dependent in the ${\mathbf Z}_p$-vector space 
${\mathbf A}\otimes{\mathbf Z}_p$ which contradicts that $b_i\,(i=1,2, \dots ,n, \dots )$ form 
a basis of ${\mathbf A}\otimes{\mathbf Z}_p$, 
where ${\mathbf Z}_p={\mathbf Z}/p{\mathbf Z}$. 
By a base change on $a_i\,(i=2,3, \dots ,n, \dots )$, 
it can be assumed that $b_2=c_{12}a_1+a_2$. 
By an inductive argument, it can be assumed that 
\[b_j=c_{1j}a_1+c_{2j}a_2+ \dots + c_{ j-1\, j} a_{j-1}+a_j\, (j=3,4, \dots ,m).\]
By a base change replacing 
$a_j$ to $a_j-c_{1j}a_1-c_{2j}a_2- \dots  -c_{j-1\, j} a_{j-1}\, (j=2,3, \dots , m)$, 
the identities $b_j=a_j\,(1\leq j\leq m)$ are obtained. 
Then the entries $c_{i j}$ of the matrix $C$ are written as 
\[ c_{11}=1\quad c_{1 j}=c_{i 1}=0\, (1< i<+\infty, 1< j< +\infty).\]
This completes the proof of Lemma~3.1. 
$\square$

\phantom{x}

The proof of the following lemma uses Lemma~3.1.

\phantom{x}

\noindent{\bf Lemma~3.2.} If the 2-complex $P$ of a locally finite group presentation 
\[GP=<x_1, x_2,  \dots , x_n,  \dots  |\, r_1, r_2,  \dots , r_m,  \dots >\]
is contractible, then 
there is a base change $x'_i\, (i=1,2, \dots )$ of the basis $x_i\, (i=1,2, \dots )$ 
of the free group $<x_1, x_2,  \dots , x_n,  \dots >$ such that the resulting group presentation 
\[GP'=<x'_1, x'_2,  \dots , x'_n,  \dots  |\, r'_1, r'_2,  \dots , r'_m,  \dots >\]
is a homology-trivial unit-group presentation.

\phantom{x}

\noindent{\bf Proof of Lemma~3.2.}
Since the 2-complex $P$ is a contractible locally finite 2-complex, 
every generator $x_i$ appears only in a finite number of the relators $r_1, r_2,  \dots , r_m,  \dots $ 
and the inclusion homomorphism 
\[<r_1, r_2,  \dots , r_m,  \dots > \to  <x_1, x_2,  \dots  x_n,  \dots >\]
on the free groups $<r_1, r_2,  \dots , r_m,  \dots >$ and $<x_1, x_2,  \dots  x_n,  \dots >$ 
induces an isomorphism on the abelianized groups 
$<r_1, r_2,  \dots , r_m,  \dots >^a$ and ${\mathbf A}=<x_1, x_2,  \dots  x_n,  \dots >^a$ 
which are free abelian groups with a base change matrix $C$ given in Lemma~3.3. 
Do Nielsen transformations on the free group $<x_1, x_2,  \dots  x_n,  \dots >$ 
induced from base changes on the free abelian group $\mathbf A$ of Lemma~3.1. 
Then the word $r_1$ is changed into $x_1$ in ${\mathbf A}$.
This base change is done by using only finitely many 
letters in $x_i\,(i=1,2, \dots ,n, \dots )$ belonging to the word $r_j$ except for re-indexing of the letters $x_i\,(i=1,2, \dots ,n, \dots )$. 
Thus, the conclusion of Lemma~3.2 is obtained. 
$\square$

\phantom{x}

The following corollary means that a contractible locally finite 2-complex may be considered as 
the 2-complex of a homology-trivial unit-group presentation.

\phantom{x}

\noindent{\bf Corollary~3.3.} There is a cellular-homotopy 
equivalence from every contractible locally finite 2-complex $P$ to the 2-complex $P'$ of 
a homology-trivial unit-group presentation $GP$ inducing a cellular-homotopy 
equivalence from the connected subcomplexes of $P$ to the subcomplexes of $GP$. 

\phantom{x}

\noindent{\bf Proof of Corollary~3.3.} 
Let $P$ be a contractible locally finite 2-complex obtained from the 1-skelton $P^1$ with 
$\pi_1(P^1,v)=<x_1, x_2,  \dots , x_n,  \dots>$ on the generators $x_i\,(i=1,2, \dots , n,  \dots ) $ 
by attaching 2-cells with attaching maps given by relators $r_j\, (j=1,2, \dots , m, \dots ) $.
Then the inclusion homomorphism 
\[<r_1, r_2,  \dots , r_m,  \dots >  \to  <x_1, x_2,  \dots  x_n,  \dots >\]
induces an isomorphism from the abelianized group 
$<r_1, r_2,  \dots , r_m,  \dots >^a$ to the abelianized group $<x_1, x_2,  \dots , x_n,  \dots >^a$.
Let 
\[g:<x_1, x_2,  \dots , x_n,  \dots > \to <x'_1, x'_2,  \dots , x'_n,  \dots >\] be a base change 
isomorphism sending the word $r_j$ to a word $r'_j$ such that $r'_j$ is equal to $x'_j$ in 
the abelianized group $<x_1, x_2,  \dots  x_n,  \dots >^a$ for all $j$. 
Let $P'$ be the 2-complex of the homology-trivial unit-group presentation 
\[GP'=<x'_1, x'_2,  \dots , x'_n,  \dots  |\, r'_1, r'_2,  \dots , r'_m,  \dots >.\]
The isomorphism $g$ induces a desired cellular homotopy equivalence $P \to  P'$.
$\square$

\phantom{x}

Conjecture~2 is obtained from Conjecture~3 as follows.

\phantom{x}

\noindent{\bf 3.4: Proof of Conjecture~3 $\mathbf \Rightarrow$ Conjecture~2.} 
By Corollary~3.3, every connected subcomplex of every contractible locally finite 2-complex 
is homotopy equivalent to a subcomplex of a homology-trivial unit-group presentation. 
This completes the proof of Conjecture~3 $\mathbf \Rightarrow$ Conjecture~2. $\square$

\phantom{x}

\noindent{\bf 4. A ribbon sphere-link and a ribbon disk-link constructed from a homology-trivial unit-group presentation}

Let $X$ be an open connected oriented smooth 4D manifold. 
A countably-infinite system of disjoint compact sets $X_i\,(i=1,2, \dots , n, \dots )$ in $X$ is
{\it discrete} if the set $\{p_i|\, i=1,2, \dots , n, \dots \}$ constructed from any one point $p_i\in X_i$ for every $i$ is a discrete set in $X$. 
An $S^2$-{\it link} in $X$ is the union $L$ of 
a discrete (finite or countably-infinite)  system of disjoint 2-spheres smoothly embedded in $X$. 
An $S^2$-link in $X$ is {\it trivial} if it bounds a discrete system of mutually disjoint 
3-balls smoothly embedded in $X$, and {\it ribbon} if it is obtained from a trivial $S^2$-link 
$O$ by surgery along a discrete system of disjoint 1-handles embedded in $X$. 
An $S^2$-{\it link} $L$ in $X$ is {\it finite} if the number of the components of $L$ is finite
Otherwise, $L$ is {\it infinite}. 
The {\it open 4D handlebody} 
\[Y^O={\mathbf R}^4\#_{i=1}^{+\infty} S^1\times S^3_i\]
denotes the connected sum of the 4-space ${\mathbf R}^4$ and a discrete system of  $S^1\times S^3_i\, (i=1,2, \dots ,n, \dots )$.
The following lemma is basic to our purpose.

\phantom{x}

\noindent{\bf Lemma~4.1.} 
For every homology-trivial unit-group presentation
\[GP=<x_1, x_2,  \dots , x_n,  \dots  |\, r_1, r_2,  \dots , r_m,  \dots >,\] 
there is a ribbon $S^2$-link $L$ with components 
$K_i\, (i=1,2, \dots ,n, \dots )$ in 
${\mathbf R}^4$ such that there is an isomorphism 
\[\pi_1(E(L),v) \to  <x_1, x_2,  \dots  x_n,  \dots >\]
sending a meridian system of $K_i\, (i=1,2, \dots ,n, \dots )$ to the relator system 
$r_i \, (i=1,2, \dots ,n, \dots )$.

\phantom{x}

\noindent{\bf Proof of Lemma~4.1.} 
In the open 4D handlebody 
$Y^O={\mathbf R}^4\#_{i=1}^{+\infty} S^1\times S^3_i$,
let $\gamma^O$ be a legged loop system  with loop system 
$k^O_i=S^1\times {\mathbf 1}_i\,(i=1,2, \dots ,n, \dots )$ representing 
a basis $x_i\,(i=1,2, \dots ,n, \dots )$ of the free group $\pi_1(Y^O,v)$.  
Let $k_j\,(j=1,2, \dots ,m, \dots )$ be a simple loop system $k_*$ in $Y^O$ 
representing the relator system $r_j\,(j=1,2, \dots ,m, \dots )$. 
By assumption of the homology-trivial unit-group presentation $GP$, 
the loop $k_j$ for every $j$ meets transversely $1\times S^3_i$ in $Y^O$ 
with intersection number $+1$ for $j=i$ and with intersection number $0$ for $j\ne i$. 
Further, the loop $k_j$ does not meet $1\times S^3_i$ except for a finite number of $i$. 
Let $X$ be the smooth open 4D manifold obtained from $Y^O$ by surgery along the loops 
$k_j\,(j=1,2, \dots ,m, \dots )$ replacing a normal $D^3$-bundle $k_j\times D^3$ of $k_j$ in 
$Y^O$ with the $D^2$-bundle $D_j\times S^2$ of $S^2$ for a disk $D_j$ with 
$\partial D_j=k_j$. 
Then the $S^2$-link $L=\cup_{j=1}^{+\infty} K_j$ with $K_j=0_j\times S^2$ is obtained in $X$. 

\phantom{x}

\noindent{\bf (4.1.1)} The open 4D manifold $X$ is contractible.

\phantom{x}

\noindent{\bf Proof of (4.1.1).} By van Kampen theorem, $X$ is simply connected because 
the loops $k_j\,(j=1,2, \dots ,m, \dots )$ normally generate the free 
fundamental group $\pi_1(Y^O,v)=<x_1,x_2, \dots ,n, \dots >$. 
Thus, if $H_q(X;{\mathbf Z})=0\, (q=2,3)$, then $X$ is contractible since 
$X$ is an open 4D manifold. By the excision isomorphism 
\[H_q(Y^O,k_*\times D^3;{\mathbf Z})\cong H_q(X, D_*\times S^2;{\mathbf Z}).\]
Hence $H_3(X, D_*\times S^2;{\mathbf Z})=0$, so that $H_3(X;{\mathbf Z})=0$. 
Since the loop system $k_*$ meets transversely $1\times S^3_i$ 
in a finite number of points in $Y^O$ with intersection number
$\mbox{Int}(k_j, 1\times S^3_i)=+1\, (j=i), 0\, (j\ne i)$,   there is 
an arc system $I_s\,(s=1,2, \dots  u)$ in the 1D manifold  system obtained from $k_*$ by cutting along the set $k_*\cap 1\times S^3_i$
such that  $I_s$ attaches  to $1\times S^3_i$ with opposite signs for all $s$ and 
the 3D  orientable manifold $Z_i$ obtained from $1\times S^3_i$ by  piping along 
$I_s\,(s=1,2, \dots  u)$ meets $k_i$ with just one point and does not meet $k_j\, (i\ne j)$.  
By the construction of $X$, the component $K_i$ of $L$ bounds 
a once-punctured 3-manifold $V_i$ of $Z_i$ in $X$ not meeting $L\setminus K_i$, for every $i$. 
This means that the inclusion homomorphism 
$H_2(D_i\times S^2;{\mathbf Z}) \to  H_2(X;{\mathbf Z})$ is the zero map for all $i$.
Thus, 
\[H_2(X;{\mathbf Z})\cong H_2(X, \cup_{i=1}^{+\infty}D_i\times S^2;{\mathbf Z})\cong 
H_2(Y^O,\cup_{i=1}^{+\infty}k_i\times D^3;{\mathbf Z})=0. \]
Thus, $X$ is a contractible open 4D manifold. 
This completes the proof of (4.1.1). 
$\square$

\phantom{x}

The proof of Free Ribbon Lemma in\cite{KA} means that the 2-sphere component $K_i$ of $L$ 
is isotopic to a ribbon $S^2$-knot in $X$ obtained from a finite trivial $S^2$-link $O_i$ split from 
$L$ by surgery along a finite disjoint 1-handle system ${\mathbf h}_i$ such that 
$\cup_{i=1}^{+\infty} O_i$ is a trivial link and ${\mathbf h}_i\,(i=1,2, \dots ,n, \dots )$ are disjoint discrete systems. 
Thus, the $S^2$-link $L$ is a ribbon $S^2$-link in $X$. 
By taking the upper-half 4-space ${\mathbf H}^4$ near the end of 
the connected summand ${\mathbf R}^4$ of $Y^O$, let $i_X:{\mathbf H}^4 \to  X$
be a smooth embedding and consider $X$ as a ${\mathbf R}^3$-connected sum
For a 4-space ${\mathbf R}^4$ in $X$, 
 the ribbon $S^2$-link $L$ in $X$  can be moved into  ${\mathbf R}^4$, 
since $L$ is obtained from a discrete trivial $S^2$-link which is movable into ${\mathbf R}^4$ 
by surgery along disjoint discrete systems ${\mathbf h}_i\,(i=1,2, \dots ,n, \dots )$ 
which are also movable into ${\mathbf R}^4$. 
By construction, there is an isomorphism from 
$\pi_1(X\setminus L,v)\cong \pi_1({\mathbf R}^4\setminus L,v)$ 
is isomorphic to the free 
fundamental group $\pi_1(Y^O,v)=<x_1,x_2, \dots ,n, \dots >$ 
sending a meridian system of $K_i\, (i=1,2, \dots ,n, \dots )$ to the relator system 
$r_i \, (i=1,2, \dots ,n, \dots )$ of $GP$. 
This completes the proof of Lemma~4.1. 
$\square$

\phantom{x}

Let 
\[{\mathbf H}^4=\{(x,y,z,t)|\, -\infty<x,y,z<+\infty, 0\leq t<\infty \}\] 
be the {\it upper-half 4-spac}e of ${\mathbf R}^4$ with boundary 
$\partial {\mathbf H}^4= \{(x,y,z,0)|\, -\infty<x,y,z<+\infty\}$ identifying 
the 3-space ${\mathbf R}^3= \{(x,y,z)|\, -\infty<x,y,z<+\infty\}$. 
Let $\alpha$ be the reflection in ${\mathbf R}^4$ sending $(x,y,z,t)$ to $(x,y,z,-t)$.
The image $\alpha(H^4)$ of the upper-half 4-space $H^4$ by $\alpha$ is given by 
the lower-half 4-space $\{(x,y,z,t)|\, 0<x,y,z<+\infty, -\infty<t\leq 0\}$. 
A {\it disk-link} $L^D$ in $H^4$ is a discrete (finite or countably-infinite) 
system of disjoint disks smoothly and properly 
embedded in ${\mathbf H}^4$. 
A (possibly infinite) disk-link $L^D$ in $H^4$ is {\it trivial} if it is obtained from 
a discrete system of disjoint disks in ${\mathbf R}^3$ by pushing up the interiors of the disks 
into the interior of ${\mathbf H}^4$. 
A  disk-link $L^D$ in $H^4$ is {\it ribbon}
if it is obtained from a disjoint discrete embedded disk system ${\mathbf D}\cup {\mathbf b}$ 
in ${\mathbf H}^4$ which is the union of 
a trivial disk-link ${\mathbf D}= \{D_i|\, i=1,2, \dots ,n, \dots \}$ in $H^4$ 
and a disjoint spanning band system ${\mathbf b}=\{b_j |\, j=1,2, \dots ,m, \dots \}$ 
in ${\mathbf R}^3$ by pushing up the interiors of the disk system 
${\mathbf D}\cup {\mathbf b}$ into the interior of ${\mathbf H}^4$. Thus, 
\[L^D=\tilde {\mathbf D}\cup \tilde{\mathbf b}\]
for a pushing up disk system $\tilde{\mathbf D}= \{\tilde D_i|\, i=1,2, \dots ,n, \dots \}$ of 
${\mathbf D}$ and a pushing up the band system
$\tilde{\mathbf b}=\{\tilde b_j |\, j=1,2, \dots ,m, \dots \}$ of ${\mathbf b}$. 
The {\it closed exterior} of a ribbon disk-link $L^D$ in ${\mathbf H}^4$  
is the 4D manifold $E(L^D)=\mbox{cl}({\mathbf H}^4\setminus N(L^D))$ for a regular neighborhood of $L^D$ in ${\mathbf H}^4$. 
Every ribbon $S^2$-link $L$ in ${\mathbf R}^4$ is isotopically deformed into an 
$\alpha$-invariant position for the reflection $\alpha$ in ${\mathbf R}^4$, so that $L$ is obtained from a  ribbon disk-link $L^D$ 
in ${\mathbf H}^4$ by doubling of ${\mathbf H}^4$ by $\alpha$ (see \cite[II]{KSS}). 
The following lemma is shown by the same method as \cite[ Lemma~3.1]{KA}.

\phantom{x}

\noindent{\bf Lemma~4.2.}
For a ribbon disk-link $L^D$ in ${\mathbf H}^4$ in a (possibly infinite) ribbon $S^2$-link $L$ in 
${\mathbf R}^4$, the inclusion $({\mathbf H}^4, L^D) \to  ({\mathbf R}^4,L)$ induces an isomorphism
\[\pi_1(E(L^D), v)  \to  \pi_1(E(L), v).\]

\phantom{x}

The following corollary is obtained from Lemmas~4.1 and 4.2.

\phantom{x}

\noindent{\bf Corollary~4.3.}
For every homology-trivial unit-group presentation
\[GP=<x_1, x_2,  \dots , x_n,  \dots  |\, r_1, r_2,  \dots , r_m,  \dots >,\] 
there is a ribbon disk-link $L^D$ with components $K^D_i\, (i=1,2, \dots ,n, \dots )$ in 
${\mathbf H}^4$ such that there is an isomorphism 
\[\pi_1(E(L^D),v) \to  <x_1, x_2,  \dots  x_n,  \dots >\]
sending a meridian system of $K^D_i\, (i=1,2, \dots ,n, \dots )$ to the relator system 
$r_i \, (i=1,2, \dots ,n, \dots )$. 

\phantom{x}

\noindent{\bf 5. A ribbon disk-link corresponding to a 1-full subcomplex of a homology-trivial unit-group presentation} 

A ribbon disk-link $L^D$ in ${\mathbf H}^4$ is {\it free} if the fundamental group 
$\pi_1(E(L^D),v)$ is a free group. 
The following lemma contains an infinite version of 
the results of \cite[Theorem~1.4, Lemma~3.2]{KA}.

\phantom{x}

\noindent{\bf Lemma~5.1.} The closed exterior $E(L^D)$ of every (possibly infinite) 
ribbon disk-link $L^D$ in ${\mathbf H}^4$ is aspherical. In particular, 
for every (possibly infinite) free ribbon disk-link $L^D$ in 
${\mathbf H}^4$ with $\pi_1(E(L^D),v)\cong <x_1,x_2, \dots ,x_n, \dots >$, there is a strong deformation retract $r: E(L^D)  \to  \omega x$ for a locally finite graph.$\omega x$ 
with $\pi_1(\omega x,v)\cong <x_1,x_2, \dots ,x_n, \dots >$. 

\phantom{x}

\noindent{\bf Proof of Lemma~5.1.}
First, show that $E(L^D)$ is always aspherical for an infinite ribbon disk-link $L^D$ in
the upper 4-space ${\mathbf H}^4$. 
Let $L^D=\tilde {\mathbf D}\cup \tilde{\mathbf b}$. 
Divide the upper 4-space ${\mathbf H}^4$ along the upper 3-space 
${\mathbf H}^3_0=\{(x,y,0,t)|\, -\infty<x,y<+\infty, 0\leq t<\infty \}$ 
into the 2-parts 
${\mathbf H}^4_+=\{(x,y,0,t)|\, -\infty<x,y<+\infty, 0\leq z, t<\infty \}$ and 
${\mathbf H}^4_-=\{(x,y,0,t)|\, -\infty<x,y<+\infty, -\infty<x\leq 0, 0\leq t<\infty \}$. 
Assume that the trivial disk-link ${\mathbf D}= \{D_i|\, i=1,2, \dots ,n, \dots \}$ 
is disjoint from ${\mathbf H}^3_0$ and splits into two disk system
${\mathbf D}_{\pm}$ 
so that ${\mathbf D}_+$ is a finite trivial disk-link in ${\mathbf H}^4_+$ 
and ${\mathbf D}_-$ is an infinite trivial disk-link in ${\mathbf H}^4_-$. 
Let $\tilde {\mathbf D}_{\pm}$ be the pushing up disk systems 
of ${\mathbf D}_{\pm}$ in ${\mathbf D}$. 
The spanning band system ${\mathbf b}$ meets ${\mathbf H}^3_0$ with a disjoint simple arc system consisting of an arc parallel to an arc attaching to the disk system ${\mathbf D}$. 
The band system ${\mathbf b}_{\pm}={\mathbf b}\cap {\mathbf H}^4_{\pm}$ 
consists of a band system ${\mathbf b}_{\pm}^1$ of bands with no end or one end in 
${\mathbf H}^3_0$ and a band system ${\mathbf b}_{\pm}^2$ of bands with both ends in 
${\mathbf H}^3_0$. 
Let $\tilde{\mathbf b}_{\pm}=\tilde{\mathbf b}_{\pm}^1\cup \tilde{\mathbf b}_{\pm}^2$ be the pushing up band systems of ${\mathbf b}_{\pm}$. 
Note that the band system ${\mathbf b}_+={\mathbf b}_+^1\cup {\mathbf b}_+^2$ 
is a finite band system. 
Let $f: S^q \to  \mbox{Int}E(L^D)$ be a map from the $q$-sphere $S^q$ for $q\geq 2$. 
By a slide of the upper 3-space ${\mathbf H}^3_0$, it can be assumed that 
the image $f(S^q)$ is in the interior of ${\mathbf H}^4_+$ and does not meet 
$L^D_+=\tilde{\mathbf D}_+\cup\tilde{\mathbf b}_+^1$ and the pushing up 
process from the band system ${\mathbf b}_+^2$ to the band system $\tilde{\mathbf b}_+^2$. 
Let $E(L^D_+)=\mbox{cl}({\mathbf H}^4_+\setminus N(L^D_+))$ 
and $E(L^D_+\cup \tilde{\mathbf b}_+^2)
=\mbox{cl}({\mathbf H}^4_+\setminus N(L^D_+\cup\tilde{\mathbf b}_+^2))$ 
for regular neighborhoods $N(L^D_+)$ and $N(L^D_+\cup \tilde{\mathbf b}_+^2)$ 
of $L^D_+$ and $L^D_+\cup \tilde{\mathbf b}_+^2$ in ${\mathbf H}^4_+$, respectively. 
Let $\bar L^D_+$ be a finite ribbon disk-link in ${\mathbf H}^4$ obtained 
from $L^D_+$ in ${\mathbf H}^4_+$ by taking a double along ${\mathbf H}^3_0$. 
Since $E(\bar L^D_+)$ is aspherical by \cite[Lemma~3.2]{KA} and there is a retraction 
$r: E(\bar L^D_+) \to  E(L^D_+)$, the inclusion $ E(L^D_+) \to  E(\bar L^D_+)$ induces 
a monomorphism $\pi_q(E(L^D_+),v) \to  \pi_q(E(\bar L^D_+),v)$. 
Thus, the map $f:S^q \to  \mbox{Int}E(L^D_+)$ defined by $f:S^q \to  \mbox{Int}E(L^D)$
extends to a map
\[f^+:D^{q+1} \to  \mbox{Int}E(L^D_+)\]
from the $(q+1)$-disk $D^{q+1}$. 
Push up the union 
$L^D_+\cup{\mathbf b}_+^2$ into $L^D_+\cup\tilde{\mathbf b}_+^2$ by a deformation keeping $L^D_+$ and $f(S^q)$ fixed. 
Since $f^+(D^{q+1})\cap {\mathbf b}_+^2=\emptyset$, this deformation only deforms the 
the part $f^+(D^{q+1})\setminus f(S^q)$ of $f^+(D^{q+1})$ isotopically. 
By construction $E(L^D_+\cup \tilde{\mathbf b}_+^2)\subset E(L^D)$, 
the image $f^+(D^{q+1})$ is in $E(L^D)$. This means that $E(L^D)$ is aspherical. 
Further, if $\pi_1(E(L^D),v)$ is isomorphic to $ <x_1,x_2, \dots ,x_n, \dots >$, then $E(L^D)$ is homotopy equivalent to $\omega x$ and there is a strong deformation retract $r: E(L^D)  \to  \omega x$.
This completes the proof of Lemma~5.1.
$\square$

\phantom{x}

For a free ribbon disk-link $L^D$ in ${\mathbf H}^4$, let 
\[Q(L^D)=E(L^D)\cup N(L^D)\] 
be a decomposition of ${\mathbf H}^4$ into the closed complement $E(L^D)$ and the normal disk-bundle $N(L^D)=L^D\times D^2$. 
Let $p_*(L^D)=\{p_i|\,i=1,2, \dots ,n, \dots  \}$ be a discrete set made by taking one point from every component of $L^D$.
The strong deformation retract $r: E(L^D)  \to  \omega x$ in Lemma~5.1
and the strong deformation retract $N(L^D) \to  p_*(L^D)\times D^2$ shrinking $L^D$ into 
$p_*(L^D)$ define a map $\rho: Q(L^D)  \to  P(L^D)$
for a connected locally finite 2-complex 
\[P(L^D)= \omega\cup p_*(L^D)\times D^2\] 
with the attaching map 
$p_*(L^D)\times \partial D^2 \to \omega x$ defined by $r$. 
The map $\rho$ is called a {\it ribbon disk-link presentation} for the 2-complex $P(L^D)$. 
For a sublink $K^D$ of $L^D$, let $N(K^D)=K^D\times D^2$ be the subbundle of the disk-bundle $N(L^D)$. The union 
\[Q(K^D;L^D)=E(L^D)\cup N(K^D)\] 
is a decomposition of the closed complement $E(L^D\setminus K^D)$ 
of the sublink $L^D\setminus K^D$  of $L^D$ in ${\mathbf H}^4$, which is a ribbon $S^2$-link 
in ${\mathbf H}^4$. The ribbon disk-link presentation 
$\rho: Q(L^D)  \to  P(L^D)$ for $P(L^D)$ sends $Q(K^D;L^D)$ to the 1-full 2-subcomplex 
\[P(K^D;L^D)=\omega\cup p_*(K^D)\times D^2\] 
of $P(L^D)$. 
Further, every 1-full 2-subcomplex of $P(L^D)$ is obtained from a sublink $K^D$ of $L^D$ in this way. 
The following theorem contains an infinite version of a ribbon disk-link $L^D$ 
of \cite[Theorem~1.3]{KA}. 

\phantom{x}

\noindent{\bf Theorem~5.2.} 
For every free ribbon disk-link $L^D$ in ${\mathbf H}^4$, 
the ribbon disk-link presentation $\rho: Q(L^D)  \to  P(L^D)$ induces a homotopy equivalence 
$Q(K^D;L^D)  \to  P(K^D;L^D)$ for every sublink $K^D$ of $L^D$ including 
$K^D=\emptyset$ and $K^D=L^D$. 
In particular, the 2-complex $P(L^D)$ is contractible. 
The 2-complex $P$ of every homology-trivial unit-group presentation $GP$
 is taken as $P=P(L^D)$ for a free ribbon disk-link $L^D$ in ${\mathbf H}^4$ 
 so that  for every 1-full subcomplex $P'$ of $P$, there is just one  sublink $K^D$ of $L^D$ 
 with  $P'=P(K^D;L^D)$. 

\phantom{x}

\noindent{\bf Proof of Theorem~5.2.} 
The homotopy equivalence of the ribbon disk-link presentation $\rho: Q(L^D)  \to  P(L^D)$
is similar to the proof of \cite[Theorem~1.3]{KA}. 

Let $GP=<x_1, x_2,  \dots , x_n,  \dots  |\, r_1, r_2,  \dots , r_m,  \dots >$ 
be a homology-trivial unit-group presentation.
By Corollary~4.3,  
there is a free ribbon disk-link $L^D$  in 
${\mathbf H}^4$ with  an isomorphism 
$\pi_1(E(L^D),v) \cong  <x_1, x_2,  \dots  x_n,  \dots >$ 
sending a meridian system of $L^D$ to the relator system $r_i \, (i=1,2, \dots ,n, \dots )$. 
The 2-complexes $P$ of $GP$ and $P(L^D)$ are both obtained from a graph $\omega x$ 
with $\pi_1(\omega x,v)=<x_1, x_2,  \dots  x_n,  \dots >$ by attaching 2-cells with attaching maps given by the relator words $r_j\,(j=1,2, \dots ,m, \dots )$. Hence the 1-full subcomplexes of $P$ 
coincide with the 1-full subcomplexes of $P(L^D)$. 
$\square$

\phantom{x}

The following corollary confirms that Conjecture~3 is true.

\phantom{x}

\noindent{\bf Corollary~5.3.} Every subcomplex of every homology-trivial unit-group presentation is aspherical. 

\phantom{x}

\noindent{\bf Proof of Corollary~5.3.} 
Let $P$ be the 2-complex  of every homology-trivial unit-group presentation, and $P'$ a 
connected subcomplex of $P$.
By Theorem~5.2, $P$ is written as $P(L^D)$ for  a free ribbon disk-link $L^D$ in 
${\mathbf H}^4$. If  $P'$ is a $1$-full subcomplex $P'$ of $P$, then $P'$ is written as 
$P(K^D;L^D)$ for a sublink $K^D$ of  $L^D$ in ${\mathbf H}^4$. 
The ribbon disk-link presentation $\rho: Q(K^D;L^D)  \to  P(K^D;L^D)$ is  homotopy equivalent  and $Q(K^D;L^D) $ is the closed exterior $E(L^D\setminus K^D)$, which is aspherical  by Lemma~5.1.  Thus, $P'$ is aspherical. 
If $P'$ of $P$ is not 1-full, then a 1-full subcomplex $P''$ of $P$ is constructed from $P'$ by adding some loops in the 1-skelton $P^1=\omega x$ to $P'$, and $P''$ is aspherical if and only if $P'$ is aspherical. Thus, $P'$ is aspherical in this case. 
This completes the proof of Corollary~5.3. 
$\square$

\phantom{x}

The proof of Theorem~1.1 is now completed as follows.

\phantom{x}

\noindent{\bf Proof of Theorem~1.1.} 
The proof of Theorem~1.1 is completed by Corollary~5.3 (a confirmation of Conjecture~3) and 
the proofs of Conjecture~3 $\mathbf \Rightarrow$ Conjecture~2 and Conjecture~2 $\mathbf \Rightarrow$ Conjecture~1. 
This completes the proof of Theorem~1.1. 
$\square$

\phantom{x}

\noindent{\bf Acknowledgments.} This work was partly supported by JSPS KAKENHI Grant
Numbers JP19H01788, JP21H00978 and MEXT Promotion of Distinctive Joint Research
Center Program JPMXP0723833165.

\phantom{x}

\end{document}